\documentclass[11pt]{amsart}
\usepackage{amsfonts}
\usepackage{latexsym}
\usepackage{amscd}
\usepackage{amsmath}
\usepackage{amssymb}
\usepackage{graphicx}
\usepackage[mathscr]{eucal}

\newtheorem{theo}{Theorem}[section]
\newtheorem{prop}{Proposition}[section]

\newtheorem{rem}{Remark}[section]
\newtheorem{co}{Corollary}[section]
\theoremstyle{definition}
\newtheorem{de}{Definition}[section]

\newcommand{\gp}{\mathbb{P}}
\newcommand{\findemo}{$\Box$ \\[1mm]}

\title{ \rm{The Abhyankar-Moh theorem for plane valuations at infinity}}
\author{C.~Galindo \and F.~Monserrat}
\curraddr{{\bf C. Galindo}: Departamento de Matem\'aticas \&
Instituto de Matem\'aticas y Aplicaciones de Castell\'on (IMAC),
Universitat Jaume I. Campus de Riu Sec, 12071 Castell\'on, Spain.
\newline {\bf F. Monserrat}: Instituto Universitario de Matem\'atica Pura y Aplicada (IUMPA), Universidad Polit\'ecnica de
Valencia, Camino de Vera s/n, 46022 Valencia, Spain.} \email{
galindo@mat.uji.es \hskip 0.3cm framonde@mat.upv.es}
\date{}
\thanks{Supported by Spain Ministry of Education
 MTM2007-64704 and Bancaixa P1-1B2009-03}

\begin{document}
\maketitle

\begin{abstract}
We introduce the class of plane valuations at infinity and prove an
analogue to the Abhyankar-Moh (semigroup) Theorem for it.
\end{abstract}

\section{Introduction}
Between 1940 and 1960, Zariski  and Abhyankar developed the theory of valuations in the context
of the theory of singularities with the aim of proving resolution
for algebraic schemes. Some important references are
\cite{zari,zar1,zar,abh1,abh2}. In the last years, there has been a
resurgence of interest in valuations within this context. Valuations of fields of quotients of
(Noetherian) regular two-dimensional  local rings $(R,m)$ centered at
$R$, which are named plane valuations, are one of the best known
classes of valuations. These valuations were classified by Zariski
%assuming that the residue field of $R$ is algebraically closed
and a refinement of that classification in terms of the dual graphs
associated with the valuations can be seen in \cite{spi} (see also
\cite{kiy}). This shows that plane valuations can be classified in
a similar way to analytically irreducible plane curve
singularities.

A particular type among these singularities is that of the
singularities at infinity of projective plane curves with only one
place at infinity. An important result that concerns this type of
curves is the Abhyankar-Moh (semigroup) Theorem \cite{pa2,pa1,am}.
It proves, under certain condition on the characteristic of the
ground field, the existence of a finite set of positive integers
satisfying certain properties that generates the
so-called semigroup at infinity and that
we call a $\delta$-sequence in $\mathbb{N}_{>0}$.
%roughly speaking, states that $\delta$-sequences span the semigroup
%at infinity of curves with only one place at infinity.
A converse result is also true: given a $\delta$-sequence in
$\mathbb{N}_{>0}$, $\Delta$, there exists a curve with only one
place at infinity whose semigroup at infinity is spanned by $\Delta$. One can
see a proof in \cite{ss} for the complex field and, for any field,
it can be deduced from \cite{cam} (see \cite{galmon} and references
therein).

The proximity between plane valuations and curve singularities and
the fact that the singularity at infinity of a curve with only one
place at infinity can be approached by approximates (see
Definition \ref{aprox})
%The fact that plane curves with only a place at infinity behave in a
%similar way to branches of plane curves
suggest the possibility of defining {\it plane valuations at infinity} (Definition \ref{preaprva}),
suitable (generalized) $\delta$-{\it sequences} and to prove for them an
analogue to the Abhyankar-Moh Theorem and its converse. We only assume that the field $k:=R/m$ is perfect. If one restricts
to valuations where $k$ is the field $\mathbb{C}$ of complex numbers and the value group is included in the set of real numbers, then plane valuations at infinity form part of the so-called valuations {\it centered} at infinity introduced in \cite{F-J-ein} to study the dynamics of polynomial mappings from $\mathbb{C}^2$ to $\mathbb{C}^2$ (see also \cite{F-J-ann}).

Coding Theory is an applied matter for which the theory of plane
curves with only one place at infinity is useful (see
\cite{far-cam} for instance). Also, good codes can be obtained if
concepts as plane valuations at infinity and attached
$\delta$-sequences are used. This was done by the authors in
\cite{galmon}; however the nature of that paper (focused on Coding
Theory) did not allow us to completely show  the extension to plane
valuations of the above mentioned results concerning plane curves
with only one place at infinity. In fact, we only proved the
converse of the Abhyankar-Moh Theorem for those types of
valuations that were useful for our purposes. This paper is
devoted to develop the mentioned suggestion in a self-contained
manner, proving an analogue to the Abhyankar-Moh Theorem for plane
valuations at infinity as they were defined in \cite{galmon} and
completing the details, concerning the definition of
$\delta$-sequence and the converse of the mentioned theorem, that
were not given in \cite{galmon}.

In Section 2, we recall the concepts of plane curve with only one
place at infinity $C$, $\delta$-sequence in $\mathbb{N}_{>0}$ and
family of approximates for curves as $C$. We also state the
Abhyankar-Moh Theorem and show the existing relation between
$\delta$-sequences in $\mathbb{N}_{>0}$ and maximal contact values
of the singularity at infinity of curves $C$ as above, fact that
will be useful in the paper. Section 3 summarizes some basic
results on plane valuations, introduces the class of plane
valuations which we are interested in, plane valuations at
infinity, and describes the value semigroup of these valuations
according their type. We normalize these semigroups to adapt them
to the concept of $\delta$-sequence for our class of valuations,
concept that is provided in Section 4. In this section, we prove
our main result, Theorem \ref{UNICO}, that we call
Abhyankar-Moh Theorem for plane valuations at infinity, whose
statement goes parallel with the classical one. We also give in Remark
\ref{r42} a nice geometrical interpretation, for divisorial valuations
at infinity, of the well-ordering of the semigroup spanned by
their associated $\delta$-sequences.

%\section{Plane valuations at infinity and $\delta$-sequences}

\section{Curves having only one place at infinity and
approximates}\label{unorr}

Along this paper, $k$ will be a perfect field and $\mathbb{P}_k^2$
(or $\mathbb{P}^2$ for short) will stand for the projective plane
over $k$.

\begin{de}
\label{de1} {\rm Let $L$ be the line at infinity in the
compactification of the affine plane to $\gp^2$. Let $C$ be a
projective absolutely irreducible curve of $\gp^2$ (i.e.,
irreducible as a curve in $\gp^2_{\overline{k}}$, $\overline{k}$
being the algebraic closure of $k$). We will say that $C$ {\it has
only one place at infinity} if the intersection $C\cap L$ is a
single point $p$ (the one at infinity) and $C$ has only one
analytic branch at $p$ (notice that the branch is defined over $k$
because $k$ is perfect)}.
\end{de}

Let $C$ be a curve with only one place at infinity. Denote by $K$
the quotient field of the local ring $\mathcal{O}_{C,p}$; the germ
of $C$ at $p$ defines a discrete valuation on $K$, that we denote by
$\nu_{C,p}$, which allows us to state the following

\begin{de}
{\rm Let $C$ be a curve with only one place at infinity given by
$p$. The {\it semigroup at infinity} of $C$ is the following
sub-semigroup of the set of non-negative integers $\mathbb{N}$
\[
S_{C,\infty}:= \left \{ - \nu_{C,p} (h) | h \in T \right \},
\]
where $T$ is the $k$-algebra $\mathcal{O}_{C}( C \setminus \{p\})$.}
\end{de}

\begin{figure}[h]
\[
\unitlength=1.00mm
\begin{picture}(80.00,30.00)(-10,3)
\thicklines \put(-5,30){\line(1,0){41}} \put(44,30){\line(1,0){16}}
\put(38,30){\circle*{0.5}} \put(40,30){\circle*{0.5}}
\put(42,30){\circle*{0.5}} \put(30,10){\line(0,1){20}}
\put(50,20){\line(0,1){10}} \put(60,0){\line(0,1){30}}
\put(10,15){\line(0,1){15}} \thinlines \put(20,30){\circle*{1}}
\put(30,30){\circle*{1}} \put(50,30){\circle*{1}}
\put(60,30){\circle*{1}}
%\put(65,30){\circle*{1}}
%\put(70,30){\circle*{1}} \put(75,30){\circle*{1}}
%\put(75,30){\circle{2}}

\put(30,20){\circle*{1}} \put(60,20){\circle*{1}}
\put(60,10){\circle*{1}} \put(10,30){\circle*{1}}
\put(30,10){\circle*{1}} \put(50,20){\circle*{1}}
\put(60,0){\circle*{1}} \put(-5,30){\circle*{1}}
\put(0,30){\circle*{1}} \put(5,30){\circle*{1}}
\put(15,30){\circle*{1}} \put(25,30){\circle*{1}}
\put(35,30){\circle*{1}} \put(45,30){\circle*{1}}
\put(55,30){\circle*{1}} \put(10,25){\circle*{1}}
\put(10,20){\circle*{1}} \put(10,15){\circle*{1}}
\put(30,25){\circle*{1}} \put(30,15){\circle*{1}}
\put(35,30){\circle*{1}} \put(-9,25){{\bf 1}=$\rho_0$}
\put(11.5,14){$\rho_1$} \put(4.5,19){$\Gamma_1$}
\put(31.5,9){$\rho_2$} \put(24.5,14){$\Gamma_2$}
\put(61.5,-1){$\rho_g$} \put(54.5,4){$\Gamma_g$} \put(9,32){$st_1$}
\put(29,32){$st_2$} \put(57.5,32){$st_g$}
%\put(72,32){$\alpha(\nu)$}
%\put(70.5,25){$\Gamma_{g+1}$}
\end{picture}
\]
\caption{The dual graph} \label{fig0}
\end{figure}
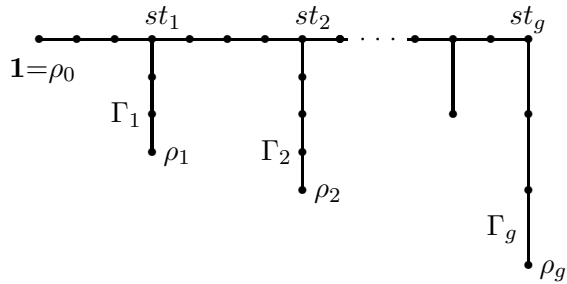

Given a curve $C$ with only one place at infinity, consider the
sequence of point blowing-ups
\begin{equation}\label{infiniteseq}
X_n\rightarrow \cdots \rightarrow X_1 \rightarrow X_0:=\gp^2
\end{equation}
that provides the minimal embedded resolution of the singularity
of $C$ at infinity. Recall that the dual graph $\Gamma$ associated
with this singularity is a tree with $n$ vertices such that each
vertex represents the strict transform in $X_n$ of one of the irreducible
exceptional divisor appearing in the above sequence of blowing-ups
and two vertices are joined by an edge whenever the corresponding
divisors intersect; additionally each vertex is labelled with the
number of blowing-ups needed to create its corresponding
exceptional divisor. The dual graph has the shape depicted in
Figure \ref{fig0}. It is the union of $g$ subgraphs $\Gamma_i$,
where $g$ is the number of characteristic pairs of the
singularity.

Let $E_{s_i}$ $(1 \leq i \leq g)$ be the exceptional divisor
obtained after blowing-up the last free point
(see the end of Section \ref{se31}) corresponding to the subgraph
$\Gamma_i$. It corresponds to the vertex $\rho_i$ in the dual
graph. An analytically irreducible germ of curve  at $p$, $\psi$,
is said to be an $i$-{\it curvette} of the germ of $C$ at $p$ if
the strict transform of $\psi$ in the surface that contains
$E_{s_i}$ is not singular and meets transversely $E_{s_i}$  and no
other exceptional curves.

From now on, for convenience, we fix homogeneous coordinates
$(X:Y:Z)$ on $\mathbb{P}^2$. $Z=0$ will be the line at infinity,
$p=(1:0:0)$ and all the singularities ``at infinity'' we will
consider will be located at $p$. Set $(x,y)$ coordinates in the
affine chart $Z\neq 0$ and $(u=y/x, v=1/x)$ coordinates around
$p$.  Note that for each polynomial $g(x,y)\in k[x,y]$, the following equality
holds:
\begin{equation}\label{equality}
g(x,y) = v^{-\deg(g)} \overline{g}(u,v),
\end{equation}
where $\overline{g}(u,v)$ is the expression of $g(x,y)$ in the local
coordinates $(u,v)$. Moreover, if $q$ is either in $k[x,y]$ or in $k[u,v]$, in this paper
by an abuse of notation, $\nu_{C,p}(q)$ will mean  $\nu_{C,p}(q')$, where $q'$ is the
element in the fraction field of $\mathcal{O}_{C,p}$ that $q$
defines.

A useful concept will be given in the following definition.

\begin{de}
\label{aprox} {\rm Let $C$ be a curve with only one place at
infinity defined by a polynomial in $k[x,y]$ of the form:
$$f(x,y)=y^m+a_1(x)y^{m-1}+\cdots+a_m(x),$$
where $m$ is the total degree of $f$, $a_i(x)\in k[x]$ ($1\leq i\leq
m$) and ${\rm char}(k)$ does not divide $m$. Let $g$ be the number
of characteristic pairs of the singularity of $C$ at $p$. A sequence of
polynomials in $k[x,y]$}, $\{q_0(x,y), q_1(x,y),  \ldots, q_g(x,y)\}$,
{\rm will be named a {\it family of approximates} for the curve $C$
whenever the following conditions happen:
\begin{description}
    \item[a] $q_0 (x,y)=x$, $q_1 (x,y)=y$,
    $\delta_0:=- \nu_{C,p} (q_0)$ and $\delta_1 :=- \nu_{C,p}
    (q_1).$
    \item[b] $q_i(x,y)$ ($1 < i \leq g$) is a monic polynomial in
    the indeterminate $y$ and $\deg(q_i)=\deg_y(q_i)=\delta_0 /d_i$, where
    $d_i=\gcd(\delta_0, \delta_1, \ldots, \delta_{i-1})$,
     being $\delta_i := - \nu_{C,p} (q_i)$.
    \item[c] If $q_i(x,y)=v^{-\deg(q_i)}\bar{q}_i(u,v)$ (see (\ref{equality})),
    then the germ of curve at $p$
    defined by $\bar{q}_i(u,v)$ (for $1<i\leq g$)
    is an $i$-curvette  (respectively, ($i-1$)-curvette) of the germ of $C$ at $p$
    when $\delta_0 - \delta_1$
    does not divide  (respectively, divides) $\delta_0$.
\end{description}
}
\end{de}

The following result is usually known as Abhyankar-Moh (semigroup)
Theorem. (See \cite{pa2,pa1,am,sathaye,pa45} as
references).

\begin{theo}
\label{21} Let $C$ be a curve having only one place at infinity and
assume that ${\rm char}(k)$ does not divide
$\gcd(-\nu_{C,p}(x),-\nu_{C,p}(y))$. Then, there exists a finite
sequence of positive integers
$\Delta=\{\delta_0,\delta_1,\ldots,\delta_g\}$ that generates the
semigroup $S_{C,\infty}$ and satisfies the following conditions:
\begin{itemize}
\item[(1)] If $d_i=\gcd (\delta_0,\delta_1,\ldots,\delta_{i-1})$,
for $1\leq i\leq g+1$, and $n_i=d_i/d_{i+1}$, $1\leq i\leq g$, then
$d_{g+1}=1$ and $n_i>1$ for $1\leq i\leq g$.

\item[(2)] For $1\leq i\leq g$, $n_i\delta_i$ belongs to the
subsemigroup of the non-negative integers,
$\mathbb{N}$, generated by $\delta_0,\delta_1,\ldots,
\delta_{i-1}$, that we usually denote $\langle
\delta_0,\delta_1,\ldots,\delta_{i-1}\rangle$.

\item[(3)] $\delta_0>\delta_1$ and $\delta_i<\delta_{i-1} n_{i-1}$
for $i=2,3,\ldots,g$.
\end{itemize}
\end{theo}

Throughout this paper, a sequence $\Delta = \{\delta_0, \delta_1,
\ldots, \delta_g\}$ of positive integers that satisfies the
conditions (1), (2) and (3) of Theorem \ref{21} will be named a
{\it $\delta$-sequence in $\mathbb{N}_{>0}$}, where
$\mathbb{N}_{>0}$ denotes the set of strictly positive integers.

\begin{rem}\label{nota1}

{\rm  The proof of Theorem \ref{21} is based on the fact that, under
the assumptions of Definition \ref{aprox} on the polynomial
$f(x,y)\in k[x,y]$ defining $C$, a specific family of approximates
for the curve (called {\it approximate roots}) can be obtained. We
point out here that, given an equation $f(x,y)=0$ for the curve $C$,
the mentioned assumptions on $f$ can be attained  by means of a change
of variables, obtaining an isomorphic curve $C'$ satisfying all the
requirements. Then, the generators of $S_{C',\infty}=S_{C,\infty}$
mentioned in the statement of Theorem \ref{21} are just the values
$\delta_0,\ldots,\delta_g$ associated with the family of
approximates given by the approximate roots of $C'$ (see Definition
\ref{aprox}). Notice also that $\delta_0$ coincides with $\deg(C')$.
}
\end{rem}

%\begin{rem}\label{nota2}
%{\rm

%If a curve $C$ as in Theorem \ref{21} satisfies the conditions
%(i), (ii) and (iii) given above, then, by applying B\'ezout's
%Theorem to $C$ and the line of infinity (whose local equation at
%$p$ is $v=0$), it happens that $\delta_0=\deg(C)$.

%}
%\end{rem}

 From Equality (\ref{equality}), Remark \ref{nota1} and the
properties that define the families of approximates, it is easy to
deduce the following result, which shows how the maximal contact
values of the singularity at infinity of a curve $C$ as in Theorem
\ref{21} (see \cite[Ch. 4]{cam}, for instance) can be expressed in
terms of a $\delta$-sequence in $\mathbb{N}_{>0}$ attached to the
semigroup $S_{C,\infty}$.

\begin{prop}\label{nuevaprop}
Let $C$ be a curve satisfying the hypotheses of Theorem \ref{21} and
let  $\{\bar{\beta}_i\}_{i=0}^g$ be the set of maximal contact
values of the germ of $C$ at $p$. Then, there exists a
$\delta$-sequence in $\mathbb{N}_{>0}$, $\{\delta_i\}_{i=0}^s$,
generating the semigroup $S_{C,\infty}$, such that if $\nu_{C,p}(u)$
does not divide (respectively, divides) $\nu_{C,p}(v)$, then $s=g$,
$\bar{\beta}_0=\delta_0-\delta_1$, $\bar{\beta}_1=\delta_0$ and
$\bar{\beta}_i=\frac{(\delta_0)^2}{d_i}-\delta_i$ for $1<i\leq g$
(respectively, $s=g+1$, $\bar{\beta}_0=\delta_0-\delta_1$, and
$\bar{\beta}_{i-1}=\frac{(\delta_0)^2}{d_i}-\delta_i$ for $1<i\leq
g+1$).
\end{prop}

The following result relates $\delta$-sequences in $\mathbb{N}_{>0}$
corresponding to two curves with only one place at infinity such
that the minimal embedded resolution of one of them is a resolution
of the other one.

\begin{co}
\label{coro} Let $C$ and $C'$ be two curves satisfying the
hypotheses of Theorem \ref{21} and such that the sequence of
morphisms (\ref{infiniteseq}) that provides a minimal embedded
resolution of the singularity of $C$ at $p$ is also an embedded
resolution of the one of $C'$. Let $\{\bar{\beta}_i\}_{i=0}^g$
(respectively, $\{\bar{\beta}'_i\}_{i=0}^{g'}\}$) be the sequence
of maximal contact values of the germ of $C$ (respectively, $C'$)
at $p$. If $\Delta=\{\delta_0,\ldots, \delta_s\}$ and
$\Delta'=\{\delta'_0,\ldots, \delta'_{s'}\}$ are
$\delta$-sequences in $\mathbb{N}_{>0}$ associated, respectively,
with the curves $C$ and $C'$ and satisfying the equalities given
in Proposition \ref{nuevaprop}, then $s'\leq s$ and
$\delta_i/\delta'_i=\bar{\beta}_0/\bar{\beta}'_0$ for $1\leq i\leq
s'-1$.
\end{co}

%\newpage

\section{Plane valuations at infinity}

\subsection{Generalities on plane valuations}
\label{se31}
To begin with, we recall the general concept of {\it valuation}.

\begin{de}
{\rm  A {\em valuation} of a field $K$ is a mapping
\[
\nu : K^* (:= K \setminus \{0\}) \rightarrow G,
\]
where $G$ is a totally ordered group, such that it satisfies
\begin{itemize}
\item $\nu (f + g) \geq \min \{\nu(f), \nu (g) \}$ and  \item $\nu
(fg) =\nu (f) + \nu (g)$,
\end{itemize}
$f, g$ being elements  in $K^*$. The subring of $K$, $R_\nu := \{ f
\in K^* | \nu(f) \geq 0\}\cup \{0\}$, is called the {\it valuation
ring} of $\nu$. $R_{\nu}$ is a local ring whose maximal ideal is
$m_\nu := \{ f \in K^* | \nu(f) > 0\} \cup \{0\}$. }
\end{de}

Given a local regular domain $(R,m)$, we will say that a valuation
$\nu$ of the quotient field of $R$ {\it is centered} at $R$ if
$R\subseteq R_{\nu}$ and $R\cap m_{\nu}=m$. The subset of $G$ given by
$S_{\nu}:=\nu(R\setminus \{0\})$ is called {\it semigroup} of the
valuation $\nu$ (relative to $R$). In the rest of the paper we will
only consider {\it plane valuations}, that is, valuations of the
quotient field of a local regular domain $(R,m)$ of dimension two
which are centered at $R$. Assume for a while that the field
$k:=R/m$ is algebraically closed. In this case, plane valuations are the algebraic version of a sequence of point
blowing-ups (see \cite{spi,kiy} for details). In fact, attached to a
plane valuation $\nu$, there is a unique sequence of point
blowing-ups
\begin{equation}
\label{uno}  \cdots \longrightarrow X_{N+1} \stackrel{\pi_{N+1}}{\longrightarrow}
X_{N} \longrightarrow \cdots \longrightarrow X_{1}
\stackrel{\pi_{1}}{\longrightarrow} X_{0}={\rm Spec} \;R,
\end{equation}
where $\pi_1$ is the blowing-up of $X_0$ centered at its closed
point $p_0$ and, for each $i\geq 1$, $\pi_{i+1}$ is the blowing-up
of $X_i$ at the unique closed point $p_i$ of the exceptional
divisor $E_i$ (obtained after the blowing-up $\pi_i$) satisfying
that $\nu$ is centered at the local ring $\mathcal{O}_{X_i,p_i}$
$(:= R_i)$. Conversely, each sequence as in (\ref{uno}) provides a
unique plane valuation. We will denote by ${\mathcal
C}_{\nu}=\{p_i\}_{i\geq 0}$ the sequence (finite or infinite) of
closed points involved in the blowing-ups of (\ref{uno}).
%this
%sequence determines the valuation $\nu$ and it is also determined
%by $\nu$.
When ${\mathcal C}_{\nu}$ is finite, $\nu$ is called the
{\it divisorial valuation} corresponding to the last exceptional
divisor  obtained in (\ref{uno}); this is so since if $\pi_{N+1}$
is the last blowing-up in the sequence (\ref{uno}) given by $\nu$,
then $\nu$ is the $m_N$-adic valuation, $m_N$ being the maximal
ideal of the ring $R_N$. Otherwise (when ${\mathcal C}_{\nu}$ is
not finite), the plane valuation $\nu$ can be regarded as the
limit of the sequence of divisorial valuations $\{\nu_i\}_{i\geq
0}$, $\nu_i$ being the divisorial valuation corresponding to
the divisor $E_i$.
%the proof is based in the fact that the ring
%$R_{\nu}$ is the direct limit of the sequence of rings
%$\{R_i\}_{i\geq 0}$.

With the above notation, let $p_i$ and $p_j$ be points in
${\mathcal C}_{\nu}=\{p_i\}_{i\geq 0}$. We will say that
$p_i$ {\it is proximate} to  $p_j$ (and it will be denoted
by $p_i\rightarrow p_j$) if $i>j$ and $p_i$ belongs to the strict
transform (by the corresponding sequence of blowing-ups given in
(\ref{uno})) of $E_{j+1}$. This binary relation among the points
of ${\mathcal C}_{\nu}$ will be called {\it proximity relation}
and it induces a binary relation ${\mathcal P}_{\nu}$ in the set of natural numbers
($i\rightarrow j$ if $p_i\rightarrow p_j$).
%it can be represented graphically by the
%so-called {\it dual graph} (see \cite{spi} for instance).
Also, the point $p_i$ is said to be {\it satellite} if there
exists $j<i-1$ such that $p_i\rightarrow p_j$ (in other words, if
$p_i$ belongs to the intersection of the strict transforms of two
exceptional divisors); otherwise, $p_i$ is said to be a {\it free}
point. It is worth pointing out that the semigroup $S_{\nu}$ of a
plane valuation depends only on the relation ${\mathcal P}_{\nu}$.
According with this relation, a plane valuation $\nu$ (with
associated sequence ${\mathcal C}_{\nu}=\{p_i\}_{i\geq 0}$)
belongs to one of the following five types (see \cite{spi} and
\cite{ga}): \vspace{3mm}

-- {\bf TYPE A}: if ${\mathcal C}_{\nu}$ is finite.

-- {\bf TYPE B}: if there exists $i_0\in {\mathbb N}_{>0}$ such
that the point $p_i$ is free for all $i>i_0$.

-- {\bf TYPE C}: if there exists $i_0\in {\mathbb N}_{>0}$ such
that $p_i\rightarrow p_{i_0}$ for all $i>i_0$.

-- {\bf TYPE D}: if there exists $i_0\in {\mathbb N}_{>0}$ such that
$p_i$ is a satellite point for all $i>i_0$ but $\nu$ is not a type C
valuation. This means that the sequence (\ref{uno})  ends with
infinitely many blowing-ups at satellite points, but they are not
ever centered at some point of the strict transforms of the same
divisor.

-- {\bf TYPE E}: if the sequence ${\mathcal C}_{\nu}$ alternates
indefinitely blocks of free and satellite points.\\

In \cite{F-J-tree}, valuations of types B, C, D and E  are  named, respectively, curve, exceptional curve, irrational and infinitely singular valuations.

\subsection{Valuations at infinity.}\label{ronaldo}

Now, we are going to define the particular type of plane valuations
that this paper deals with: {\it plane valuations at infinity}.
%Although this type of valuations has already been introduced in
%\cite{galmon}, we recall it here for convenience of the reader.
In the rest of the paper {\it we will not assume  that the ground
field $k$ is algebraically closed}. It is important to notice
that, in spite of this, the  procedure and concepts above
explained will work similarly because, due to the especial nature
of the valuations that we will consider, the centers of the
associated blowing-ups will be defined over $k$. We recall that
the field $k$ is perfect.

We start by stating the concept of general element of a
divisorial valuation.

\begin{de}
{\rm Let $\nu$ be a divisorial valuation. An element $f$ in the
maximal ideal of $R$ is named to be a {\it general element of $\nu$}
if the germ of curve given by $f$ is analytically irreducible, its
strict transform in the last variety obtained by the sequence
(\ref{uno}) attached to $\nu$, $X_{N+1}$, is smooth and meets
$E_{N+1}$ transversely at a non-singular point of the exceptional
divisor of the sequence (\ref{uno}).}
\end{de}

\begin{rem}\label{nota}
{\rm General elements are useful to compute plane divisorial
valuations. Indeed, if $f \in R$, then
$$
\nu (f) = \mbox{ min $\{(f,g) | g $ is a general element of $\nu
\},$ }
$$
where $(f,g)$ stands for the intersection multiplicity of the germs
of curve given by $f$ and $g$. The above minimum is attained when
the strict transforms of the germs defined by  $f$ and $g$ do not meet $E_{N+1}$ at the same
point.

}
\end{rem}

\begin{rem}\label{nota32}
{\rm Assume that $\{u,v\}$ is a regular system of parameters of the
ring $R$. Attached to a divisorial valuation $\nu$ and with respect
to  the system $\{u,v\}$, there exists an expression, named
Hamburger-Noether expansion of $ \nu$, that provides regular system
of parameters for the rings $R_i$ determined by $\nu$.  From this
Hamburger-Noether expansion, it can be obtained  elements $u(t,s),
v(t,s)$ in the formal power series ring in two indeterminates
$k[[t,s]]$ such that, for any $f \in R$, $\nu(f) =$ ord$_t f
\left(u(t,s), v(t,s) \right)$ (see, \cite[Sect. 3.1]{belga}).
Moreover, if $\bar{q}_i:= \bar{q}_i(u,v)$, $1 \leq i \leq g$, is a
$i$-curvette of any general element of $\nu$, then $\bar{q}_i \left
((u(t,s), v(t,s)\right) = \lambda_i s^{a_i} t^{\nu(\bar{q}_i)} +
r_i$, where $0 \neq \lambda_i \in k$, $a_i \geq 0$ and $r_i$ is an
element in $k[[t,s]]$ with exponents of $t$ larger than
$\nu(\bar{q}_i)$. Even more, for general elements of $\nu$,
$\bar{q}_{g+1} := \bar{q}_{g+1}(u,v)$, it happens that
$\bar{q}_{g+1} \left((u(t,s), v(t,s)\right) = p(s) s^a
t^{\nu(\bar{q}_{g+1})} + r_{g+1}$, $a \geq 0$ and, here, $p(s) \in
k[s] \setminus k$, $p(0) \neq 0$ and the exponents of $t$ in
$r_{g+1}$ are always larger than  $\nu(\bar{q}_{g+1})$. The
mentioned non-negative integer values $a_i$ and $a$ vanish when the
defining divisor of $\nu$ is given by a  free point. As a reference
one can see the proof of \cite[Th. 1]{belga}. There, it is
considered the free case; the behavior in the non-free case was
pointed out to the first author by A. Nu\~{n}ez when, together with
F. Delgado, they were carrying out a joint work.

%throughout the
%preparation of a joint work with her and F. Delgado.
 }
\end{rem}

Let $p$ be a closed point of $\gp^2$ on the line of infinity and
assume, from now on, that $R={\mathcal O}_{{\mathbb P}^2,p}$ and $K$
is the quotient field of $R$.

\begin{de}
{\rm A {\it plane divisorial valuation at infinity} is a plane
divisorial valuation of $K$ centered at $R$ that admits, as a
general element, an element in $R$ providing the germ at $p$ of some
curve with only one place at infinity ($p$ being its point at
infinity). }
\end{de}

\begin{de}
\label{preaprva} {\rm A plane valuation $\nu$ of $K$ centered at $R$
is said to be {\it  at infinity} whenever it is a limit of plane
divisorial valuations at infinity. More explicitly, $\nu$ will be at
infinity if there exists a sequence of divisorial valuations at
infinity $\{\nu_i\}_{i=1}^{\infty}$ such that ${\mathcal
C}_{\nu_i}\subseteq {\mathcal C}_{\nu_{i+1}}$ for all $i \in
{\mathbb N}_{>0}$ and ${\mathcal C}_{\nu}=\bigcup_{i\geq 1}
{\mathcal C}_{\nu_i}$.
% set
%$\{\nu_i\}_{i=1}^\infty$ the set of plane divisorial valuations,
%corresponding to divisors $E_i$, that appear in the sequence
%(\ref{uno}) given by $\nu$. $\nu$ will be at infinity if, for any
%index $i_0$, there is some $i > i_0$ such that $\nu_i$ is a plane
%divisorial valuation at infinity.
}
\end{de}

There exist plane valuations at infinity of all types above
described. The concept of valuation at infinity of type A is
equivalent to the one of plane divisorial valuation at infinity;
such a valuation is obtained whenever the sequence
$\{\nu_i\}_{i=1}^{\infty}$ given in the above definition satisfies
that $\nu_i=\nu_{i+1}$ for every index larger than or equal to a
fixed index $i_0\in {\mathbb N}_{>0}$ (in fact, it can be taken
constant for all $i$). It is obtained a valuation at infinity of
type B if there exists $i_0\in {\mathbb N}_{>0}$ such that
${\mathcal C}_{\nu_{i_0}}$ is the set of centers of the
blowing-ups corresponding with the minimal embedded resolution of
the germ at $p$ of a curve having only one place at infinity and,
for all $i\geq i_0$, the strict transform of this germ meets
transversely the exceptional divisor associated with $\nu_i$.
Explicit constructions of plane valuations at infinity of types C,
D, and E are described in \cite{galmon}.

The semigroup $S_{\nu}$ of a {\it valuation of type A} is, up to
equivalence of valuations, minimally  generated by the  maximal
contact values of the germ at $p$ defined by a general element of
$\nu$, $\{\bar{\beta}_i\}_{i=0}^g$; $g$ is the number  of characteristic pairs of
a general element of $\nu$ (see Sections 6 and 8 of
\cite{spi}). These maximal contact values can be computed by
evaluating by $\nu$ the $g+1$ first elements of any {\it
generating sequence} $\{\mathfrak{q}_i\}_{i=0}^{g+1}$ of $\nu$
\cite{spi}. Each element $\mathfrak{q}_i\in R$ ($0\leq i\leq g+1$)
defines an analytically irreducible germ of curve $\psi_i$ at $p$,
$\mathfrak{q}_{g+1}$ is a general element of $\nu$, and, for each
$i\leq g$, $\psi_i$ is an $i$-curvette of $\psi_{g+1}$. Notice
that, due to Remark \ref{nota}, $S_{\nu}$ coincides with the
semigroup of values of $\psi_{g+1}$. Usually, the elements in
the set $\{\nu(\mathfrak{q}_i) \}_{i=0}^{g+1}$ are named the
maximal contact values of the valuation $\nu$.

{\it Plane valuations at infinity $\nu$ of type} B correspond to
valuations at infinity in the {\it Case 4.1.a} of \cite[Sect.
9]{spi}. Here we can identify the group $G$ with ${\mathbb Z}^2$
(with the lexicographical ordering) and, with this identification,
$S_{\nu}$ is  minimally generated by the set
$\{\bar{\beta}_0:=(0,\bar{\beta}'_0),
\bar{\beta}_1:=(0,\bar{\beta}'_1),\ldots,
\bar{\beta}_g:=(0,\bar{\beta}'_g), \bar{\beta}_{g+1}:=(1,a)\}$,
with $a\in {\mathbb Z}$ (notice that $(1,a)$ is the minimum
element of $S_{\nu}$ with non-zero first coordinate), where
$\{\bar{\beta}'_i\}_{i=0}^{g+1}$ are the maximal contact values
associated with the divisorial valuation corresponding to the last
satellite point of ${\mathcal C}_{\nu}$ \cite[1.10]{d-g-n}. Since
we have freedom to choose the value of $a$, we will take, for
convention, $a=0$.

{\it Plane valuations at infinity $\nu$ of type} C correspond to
valuations at infinity in the {\it Case 3} of \cite[Sect. 9]{spi}.
The semigroup $S_{\nu}$ can be identified with a sub-semigroup of
${\mathbb Z}^2$ (with the lexicographical ordering) whose minimal
set of generators $\{\bar{\beta}_i\}_{i=0}^{g}$ can be
represented, as in the previous case, in a particular form
and with some degree of freedom (see \cite[Sect. 9]{spi} and also
\cite[1.10]{d-g-n}). We fix, for convention, the computation of
these generators in the way that we describe below. Let $i_0$ be
the maximum among the integers $j$ corresponding to points $p_j\in
{\mathcal C}_{\nu}$ that admit more than one point proximate to it
(that is, $p_{i_0}$ has infinitely many proximate points in
${\mathcal C}_{\nu}$). Since a valuation is determined by its
minimum values at the maximal ideals $m_j$ of the local rings
$R_j={\mathcal O}_{X_j,p_j}$ (denoted by $\nu(m_j)$), we set
$\nu(m_j):=(0,1)$ for all $j>i_0$, $\nu(m_{i_0}):=(1,0)$ and
$\nu(m_j):=\sum_{k\rightarrow j} \nu(m_k)$ for $j<i_0$. Now, given
a generating sequence $\{\mathfrak{q}_i\}_{i=0}^{g+1}$ of the
divisorial valuation associated to whichever divisor $E_{j}$ with
$j>i_0+2$, we have that
$\bar{\beta}_i=\nu(\mathfrak{q}_i)=\sum_{p_k\in {\mathcal
C}_{\nu}} e_k(\mathfrak{q}_i) \nu(m_k) $, $0\leq i\leq g$, where
$e_k(\mathfrak{q}_i)$ denotes the multiplicity of the strict
transform of the germ given by $\mathfrak{q}_i$ at the point
$p_k$ and the product $e_k(\mathfrak{q}_i) \nu(m_k)$ is the usual
one between elements in $\mathbb{Z}$ and $\mathbb{Z}^2$.

{\it Plane valuations at infinity $\nu$ of type} D correspond to
valuations at infinity in the {\it Case 2} of \cite[Sect. 9]{spi}.
The semigroup $S_{\nu}$ can be identified with a sub-semigroup of ${\mathbb R}$. If
$\{\nu_i\}_{i=1}^{\infty}$ is a sequence of divisorial valuations at
infinity defining $\nu$, then there exists $i_0\in {\mathbb N}$ such
that, for all $i\geq i_0$, the sets of minimal generators
$\{\bar{\beta}_j^i\}_{j=0}^g$ of $S_{\nu_i}$ have the same
cardinality $g+1$. For our future development, it is convenient to normalize the valuations $\nu_i$ by
dividing by $\nu_i(v)-\bar{\beta}_0^i$. If $g \geq 2$, set
\begin{equation}
\label{tipoD}
\frac{\bar{\beta}_j^i}{\nu_i(v)-\bar{\beta}_0^i}=\frac{\bar{\beta}_j^k}{\nu_k(v)-\bar{\beta}_0^k}=:\bar{\beta_j}\in
{\mathbb Q}
\end{equation}
for $j<g$ and $i, k > i_0$, and $\bar{\beta}_g:=\lim_{i\rightarrow
\infty} \frac{\bar{\beta}_{g}^i}{\nu_i(v)-\bar{\beta}_0^i}$ (that
belongs to ${\mathbb R}\setminus {\mathbb Q}$). Otherwise, we must
set $\bar{\beta}_j:=\lim_{i\rightarrow \infty}
\frac{\bar{\beta}_{j}^i}{\nu_i(v)-\bar{\beta}_j^i} \in {\mathbb
R}\setminus {\mathbb Q}$ $(j=0,1)$. Then the semigroup $S_{\nu}$ is generated
by $\{\bar{\beta}_i\}_{i=0}^g$.

Notice that the value semigroup of this type of valuations is
spanned by a set of finitely many positive rational numbers plus a positive
real non-rational number, fact that also happens when $g=1$ because
one can consider the isomorphic normalized semigroup generated by 1
and $ \bar{\beta}_{1}/\bar{\beta}_{0}$.

{\it Plane valuations at infinity of type} E correspond to
valuations at infinity in {\it Case 1} of \cite[Sect. 9]{spi}. The
semigroup $S_{\nu}$ can be identified with a sub-semigroup of
${\mathbb Q}$. If $\{\nu_i\}_{i=1}^{\infty}$ is a sequence of
divisorial valuations at infinity defining $\nu$, then $S_{\nu}$
is generated by an infinite set $\Gamma := \{\bar{\beta}_0,
\bar{\beta}_1,\bar{\beta}_2,\ldots\}$. We can pick a suitable
infinite sub-sequence  $\{\nu_{i_j}\}_{j=1}^{\infty}$ of
$\{\nu_i\}_{i=1}^{\infty}$ that also converges to $\nu$ and
satisfies that for each $j\geq 1$, there exists a subset of $\Gamma$,
$\{\bar{\beta}_0,\bar{\beta}_1,\ldots,\bar{\beta}_g\}$,
with $g=g(j)$ and $g(j)>g(j')$ if $j>j'$,
such that $\bar{\beta}_l=\frac{\bar{\beta}'_l}{\nu_{i_j}(v) -
\bar{\beta}'_0}$, $0\leq l\leq g$, $\{\bar{\beta}'_l\}_{l=0}^g$
being the minimal set of generators of $S_{\nu_{i_j}}$ (notice
that, for different indices $j$, the mentioned quotients must
coincide for the common indices $l$ and that this fact allows us
to normalize the valuations $\nu_{i_j}$ by dividing by
$\nu_{i_j}(v) - \bar{\beta}'_0$). For simplicity's sake we will
only consider sequences as $\{\nu_{i_j}\}_{j=1}^{\infty}$ to
define this type of valuations.

\begin{rem}
{\rm In the next section, we will define the concept the
$\delta$-sequence, which is suitable to treat the generation of
semigroups at infinity of valuations at infinity. For this purpose,
we will need, in some cases, to normalize $\delta$-sequences in
${\mathbb N}_{>0}$ and also their corresponding maximal contact
values. This fact explains why we have normalized the above sets
of maximal contact values in the described manner. }
\end{rem}

\section{Semigroup at infinity and $\delta$-sequences.}

Firstly, we introduce the concept of semigroup at infinity of a
plane valuation at infinity.

\begin{de}
{\rm Let $\nu:K^* \rightarrow G$ be a plane valuation at infinity.
The {\it semigroup at infinity} of $\nu$ is defined to be the
following sub-semigroup of $G$:
$$S_{\nu,\infty}:=\{-\nu(f)\mid f\in k[x,y] \setminus\{0\}\}.$$
}
\end{de}

Our objective is to prove an analogue to the Abhyankar-Moh Theorem
for semigroups at infinity as above. The natural first step is to
consider, in the new stage, the analogue of the $\delta$-sequences
in $\mathbb{N}_{>0}$ that appear in the Abhyankar-Moh Theorem (we
will call them simply {\it $\delta$-sequences}).
%Notice that the term
%$\delta$-sequence is used in the literature for sequences satisfying
%the conditions of the Abhyankar-Moh Theorem (here named
%$\delta$-sequences in $\mathbb{N}_{>0}$), however, for simplicity,
%we will use here this term for the corresponding sequences attached
%to plane valuations.

As a previous concept, a {\it normalized $\delta$-sequence in
$\mathbb{N}_{>0}$} will be an ordered finite set of rational
numbers $\overline{\Delta} = \{
\overline{\delta}_0,\overline{\delta}_1,\ldots,\overline{\delta}_g
\}$ such that there is a $\delta$-sequence in ${\mathbb N}_{>0}$,
$\Delta = \{ \delta_0,\delta_1,\ldots,\delta_g \}$, satisfying
$\overline{\delta}_i = \delta_i /\delta_1$ for $0 \leq i \leq g$.
Also, for any $\delta$-sequence in $\mathbb{N}_{>0}$, $\Delta = \{
\delta_0,\delta_1,\ldots,\delta_g \}$, we define the following
associated positive integers: if $\delta_0-\delta_1$ does not
divide $\delta_0$, then
\begin{equation*}
e_0:=\delta_0-\delta_1, \;\;\; e_i:=d_{i+1}
\end{equation*}
\begin{equation*}
m_0:=\delta_0, \;\;\; m_i:=n_i\delta_i-\delta_{i+1}
\end{equation*}
for $1\leq i\leq g-1$.  Otherwise,
\begin{equation*}
e_0:=d_2=\delta_0-\delta_1, \;\;\; e_i:=d_{i+2}
\end{equation*}
\begin{equation*}
m_0:=\delta_0+n_1\delta_1-\delta_2, \;\;\;
m_i:=n_{i+1}\delta_{i+1}-\delta_{i+2}
\end{equation*}
for $1\leq i\leq g-2$.

Now, we are going to give the definition of $\delta$-sequence for
the different types of plane valuations at infinity. This
definition will be followed by some examples and an explanatory
remark.
\begin{de}
{ \label{buena} A {\it $\delta$-sequence of} {\bf TYPE A}
(respectively, {\bf B}, {\bf C}, {\bf
D}, {\bf E}) is a sequence $\Delta = \{
\delta_0,\delta_1,\ldots,\delta_i, \ldots \}$ of elements in
${\mathbb Z}$ (respectively, ${\mathbb Z}^2$,
${\mathbb Z}^2$, ${\mathbb R}$,
${\mathbb Q}$) such that

\begin{description}
\item[{\bf TYPE A}] $\Delta = \{\delta_0,\delta_1,\ldots,\delta_g,
\delta_{g+1}\} \subset \mathbb{Z}$ is finite, the elements of the
set $\{\delta_0,\ldots,\delta_g\} $ satisfy the conditions
(1), (2) and (3) of the Abhyankar-Moh Theorem and $\delta_{g+1}\leq
n_g \delta_g$.

\item[{\bf TYPE B}] There exists a $\delta$-sequence in
$\mathbb{N}_{>0}$, $\Delta^* = \{
\delta_0^*,\delta_1^*,\ldots,\delta_g^* \}$, such that
$\Delta=\{(0,\delta_0^*),(0,\delta_1^*),\ldots,(0,\delta^*_g),
(-1,(\delta_{0}^*)^2)\}$.
%where $\delta_{g+1}^* =
%n_{g}^*\delta_{g}^*$.

 \item[{\bf TYPE C}] $\Delta = \{
\delta_0,\delta_1,\ldots,\delta_g \}
 \subset \mathbb{Z}^2$ is finite, $g \geq 2$ (respectively, $\geq 3$) and there exists a
 $\delta$-sequence in $\mathbb{N}_{>0}$,
 $\Delta^* = \{ \delta_0^*,\delta_1^*,\ldots,\delta_g^* \}$, such
 that $\delta_0^*-\delta_1^*$ does not divide (respectively, divides)
 $\delta^*_0$ and
\[
\delta_i = \frac{\delta_i^*}{A a_t +B} (A,B) \;\;\; (0\leq i \leq
g-1)\;\; \mbox{ and}
\]
\[
\delta_g = \frac{\delta_g^*+A' a_t +B' }{A a_t +B} (A,B) -(A',B'),
\]
where $\langle a_1;a_2, \ldots,a_t\rangle$, $a_t\geq 2$,  is the
continued
 fraction expansion of the quotient $m_{g-1}/ e_{g-1}$
 (respectively, $m_{g-2}/ e_{g-2}$) given by
 $\Delta^*$ and, considering the finite recurrence
 relation $\underline{y}_i = a_{t-i}\underline{y}_{i-1} +
 \underline{y}_{i-2}$, $\underline{y}_{-1}=(0,1)$,
 $\underline{y}_{0}=(1,0)$, then $(A,B) := \underline{y}_{t-2}$ and $(A',B') : =
 \underline{y}_{t-3}$. We complete this definition by adding that $\Delta = \{
 \delta_0,\delta_1\}$ (respectively, $\Delta = \{
 \delta_0,\delta_1,\delta_2\}$) is a $\delta$-sequence of type C
 whenever $\delta_0 = \underline{y}_{t-1}$ and $\delta_0 - \delta_1 =
 \underline{y}_{t-2}$ (respectively, $\delta_0=j\underline{y}_{t-2}$,
 $\delta_0-\delta_1=\underline{y}_{t-2}$ and $\delta_0+n_1\delta_1-\delta_2=\underline{y}_{t-1}$) for the above recurrence attached to a
 $\delta$-sequence in $\mathbb{N}_{>0}$,
 $\Delta^* = \{ \delta_0^*,\delta_1^*\}$ (respectively, $\Delta^* = \{
 \delta_0^*,\delta_1^*,\delta_2^*\}$, such that $j:=\delta_0^*/(\delta_0^*-\delta_1^*)\in \mathbb{N}_{\geq 0}$ and $n_1:=\delta_0^*/\gcd(\delta_0^*,
 \delta_1^*)$).

\item[{\bf TYPE D}] $\Delta = \{ \delta_0,\delta_1,\ldots,\delta_g
\} \subset \mathbb{R}$ is
        finite, $g \geq 2$, $\delta_i$ is a positive rational number for $0 \leq i \leq
        g-1$, $\delta_g $  is non-rational,
        and there exists a sequence $$ \left \{ \overline{\Delta}_j =
        \{ \delta_0^j,\delta_1^j,\ldots,\delta_g^j \} \right\}_{j
        \geq1}$$
        of normalized $\delta$-sequences in $\mathbb{N}_{>0}$ such
        that $\delta_i^j = \delta_i$ for $0 \leq i \leq
        g-1$ and any $j$ and  $\delta_g = \lim_{j \rightarrow \infty}
        \delta_g^j$. We complete this definition by adding that
        $\Delta = \{\tau,1\}$, $\tau>1$ being a non-rational number,
        is also a $\delta$-sequence of type D.

\item[{\bf TYPE E}] $\Delta = \{\delta_0,\delta_1,\ldots,\delta_i,
\ldots\} \subset \mathbb{Q}$ is infinite and any
        ordered subset $\Delta_j = \{\delta_0,\delta_1,\ldots,\delta_j\}$
        is a normalized $\delta$-sequence in $\mathbb{N}_{>0}$.
\end{description}
}
\end{de}

 \noindent{\bf Examples.} Next, we show  examples of
 $\delta$-sequences of types from A to E: $\{18,12,33,4,-5\}$ is of type
 A, $\{(0,18),(0,12),(0,33),(0,4), (-1,18^2)\}$ of type B, $\{(6,6),
 (4,4), (11,11), (1,2)\}$ of type C, $\{3/2,1,33/12,4/12, (33+ 14 \sqrt{2})/6(7 + 3\sqrt{2}) \}$ of type
 D and the first terms of a $\delta$-sequence of type E are $\{3/2,1,33/12,1/3,15/4,
 \ldots\}$.

 \begin{rem}
{\rm The proximity relation attached to the minimal embedded
resolution of the singularity at infinity of  a curve with only
one place at infinity that satisfies the Abhyankar-Moh Theorem can
be recovered from the finite sequence of positive integers
satisfying (1), (2) and (3) of that theorem by using the formulae
before Definition \ref{buena}. To do it, one needs to consider the
continued fractions of quotients of the type $m_l/e_l$, $m_l$ and
$e_l$ being the values defined before Definition \ref{buena} (see
Section 2 in \cite{galmon} and references therein).

$\delta$-sequences in Definition \ref{buena} are defined in such a
way that the same property happens for any type of valuation,
although we need, for that purpose, to use an extended version of
the Euclidian Algorithm that can also involve values either in
$\mathbb{Z}^2$ or in $\mathbb{R}$. We must add that in  the case
of a $\delta$-sequence of type A, $\{\delta_0,\delta_1,\ldots,
\delta_g,\delta_{g+1}\}$, the quotients $m_l/e_l$ (for $0\leq
l\leq g-1$ if $\delta_0-\delta_1$ does not divide $\delta_0$, and
for $0\leq l\leq g-2$ otherwise) determine the proximity relation
of the points of ${\mathcal C}_{\nu}=\{p_i\}_{i=0}^n$
corresponding to the minimal embedded resolution of the germ given
by a general element of $\nu$; that is, the proximity relation
among those points $p_i$ such that $i\leq i_0$, $i_0$ being the
maximal value such that $p_{i_0}$ is a satellite point. To recover
the dual graph of $\nu$, in addition to the above information, we
need to know the number $\mathfrak{f}= n-i_0$ of last free points.
This number is given by the element $\delta_{g+1}$ of the
$\delta$-sequence; in fact,
$\mathfrak{f}=n_g\delta_g-\delta_{g+1}$ (see the proof of the
forthcoming Theorem \ref{UNICO}).

Since any plane valuation $\nu$ can be regarded as a limit of
divisorial ones and since, in this last case, the continued
fractions of the quotients $m_l/e_l$ (together with
$\mathfrak{f}$) provide the corresponding proximity relations, the
definition of $\delta$-sequence for the remaining types of
valuations is made to extend at infinity the previous behavior.
So, when one can obtain  infinitely many quotients $m_l/e_l$, the
corresponding valuation is of type E. Otherwise, only finitely
many quotients $m_l/e_l$ appear and the limit of divisorial
valuations can be seen in the last one $m_r/e_r$. In these cases,
the quotients $m_l/e_l$ of the divisorial valuations $\nu_i$
converging $\nu$,  $i \gg 0$, are the same that those for $\nu$
whenever $l < r$ and the  continued fractions $\langle a_1;a_2,
\ldots,a_s\rangle$ corresponding to the last pair $(m_r,e_r)$ of the
valuations $\nu_i$ must be taken ``at infinity". This can be done either increasing the value
$s$ indefinitely or increasing $a_s$ indefinitely. In the first
case $m_r/e_r$ converges to a real non-rational number and we get
type D valuations. In the second case, we obtain type C valuations
when $s \neq 1$ and type B ones, otherwise. Note that to take $a_1$
at infinity is the same thing that to do so with the value
$\mathfrak{f}$.

Next, we make explicit for the above examples what we have said in
the previous paragraphs. The proximity relation provided by the type
A valuation is given by the following pairs $(m_l,e_l)$: $(21,6),
(62,3), (2,1)$ with attached continued fractions $\langle 3;2
\rangle$, $\langle 20;1,2 \rangle$ and $\langle 2 \rangle$; in
addition, the equality $n_g\delta_g-\delta_{g+1}=17$ indicates that
the last 17 points of ${\mathcal C}_{\nu}$ are free. For the type B
valuation, we get pairs $((0,21),(0,6)), ((0,63), (0,3))$ and
$((1,0),(0,1))$, with associated continued fractions $\langle 3;2
\rangle$, $\langle 20;1,2 \rangle$ and $\langle \infty \rangle$;
this last one corresponds to blowing-up at infinitely many free
points. In the type C case, we obtain $((7,7),(2,2))$ and $((21,20),
(1,1))$, with continued fraction $\langle 3;2 \rangle$ and  $\langle
20;1, \infty \rangle$ (indeed, for the last pair, performing the
(generalized) Euclidian Algorithm we have
$(21,0)=\mathbf{20}(1,1)+(1,0)$; $(1,1)=\mathbf{1}(1,0)+(0,1)$ and
$(1,0) = \mathbf{\infty} (0,1)$). Notice that here, $\infty$
indicates the existence  of infinitely many satellite point
blowing-ups. In case D, the pairs are $(21/12,1/2)$, $(62/12,3/12)$
and $((9+ 4 \sqrt{2})/6(7 + 3\sqrt{2}),1/6)$, and the associated
continued fractions are $\langle 3;2 \rangle$, $\langle 20;1,2
\rangle$ and $\langle 1;3,2, \sqrt{2} \rangle$. Finally, the first
pairs of the type E valuation reproduce the behavior of the one of
type A.}
\end{rem}
%Moreover, it is clear that, if  the union of the
%families of approximates of all the curves $C_{\Delta_j}$ is a
%generating sequence of $\nu_{\mathbb{G}}$.

%Assume that $\nu_{\mathbb{G}}$ is neither of type A nor of type B.
%Let $\{\bar{\beta}_k\}_{k=0}^r$ ($0\leq r\leq \infty$) be the%
%minimal system of generators of the semigroup
%$S_{\nu_{\mathbb{G}}}$. It is deduced in \cite[Sect. 4]{galmon} that
%there exists a $\delta$-sequence
%$\Delta=\{\delta_0,\delta_1,\ldots\}$ (of the same type than the one
%of the valuation) such that

%\begin{lem}
%Let $\nu$ be a plane valuation at infinity having admitting a
%generating sequence $\{\overline{Q}_i(u,v)\}_{i=0}^r$, $0\leq r\leq
%\infty$, such that $\overline{Q}_i\in k[u,v]$ for all $i\leq r$.
%Then, $S_{\nu,\infty}$ is generated by $\{-\nu(Q_i(x,y))\}_{i=0}^r$,
%where $Q_i(x,y)=v^{-\deg(Q_i)}\overline{Q}_i(u,v)$ (that is,
%$Q_i(x,y)$ defines the curve given by $\overline{Q}_i(u,v)$ in the
%chart $Z\not=0$).
%\end{lem}

To end this paper, we state and prove the announced analogue of
the Abhyankar-Moh Theorem for plane valuations at infinity.

\begin{theo}
\label{UNICO} Let $\nu$ be a plane valuation at infinity under the
assumption that it can be defined by a sequence of plane divisorial
valuations at infinity $\{\nu_i\}_{i=1}^{\infty}$ such that the characteristic of the field $k$, ${\rm
char}(k)$, does not divide $\gcd(-\nu_i(x),-\nu_i(y))$ for all $i\in
{\mathbb N}_{>0}$. Then, there exists a $\delta$-sequence $\Delta$
of the same type as $\nu$ such that the semigroup at infinity
$S_{\nu,\infty}$ is generated by $\Delta$.
\end{theo}

\noindent {\it Proof}. We will prove the result by showing it for
each type of valuation. First, assume that {\it $\nu$ is a valuation
of type} A and consider a curve $C$ having only one place at
infinity whose germ at $p$ defines a general element of $\nu$. Since
$\nu(x)=\nu_{C,p}(x)$ and $\nu(y)=\nu_{C,p}(y)$, applying Theorem
\ref{21}, the existence of a $\delta$-sequence
$\Gamma=\{\delta_0,\delta_1,\ldots,\delta_g\}$ in $\mathbb{N}_{>0}$
generating the semigroup $S_{C,\infty}$ is proved. If ${q}_{g+1}(x,y)=0$
is the affine equation (in the chart
$Z\not=0$) of the curve $C$, without loss of generality, we can assume
that the polynomial $q_{g+1}$ satisfies the
requirements of Definition \ref{aprox}. Let $\{q_i(x,y)\}_{i=0}^g$
be a family of approximates for $C$ such that
$\delta_i=-\nu_{C,p}(q_i)$, $0\leq i\leq g$.  As usual, we will
identify polynomials in $k[x,y]$ with their images in ${\mathcal
O}_{C}(C \setminus \{ p \})$.

Notice that $\nu(\bar{q}_i(u,v))=\nu_{C,p}(\bar{q}_i(u,v))$ when
$0\leq i\leq g$ and, as a consequence, the inclusion
$\Gamma\subseteq S_{\nu,\infty}$ holds. Let $i_0$ be the largest
index such that the point $p_{i_0}\in {\mathcal
C}_{\nu}=\{p_i\}_{i=0}^n$ is satellite, then
$$-\nu({q}_{g+1}(x,y))=\delta_0\nu(v)-\nu(\bar{q}_{g+1}(u,v))=\delta_0^2-n_g\nu(\bar{q}_g(u,v))-(n- i_0)=$$ $$\delta_0^2-n_g\left(\frac{\delta_0^2}{n_g}+\nu(q_g(x,y))\right)-(n- i_0)=n_g\delta_g-(n-i_0).$$
If we set $\delta_{g+1}:=n_g\delta_g-(n-i_0)$, it is clear that
$\Delta := \{\delta_0,\delta_1,\ldots,\delta_g,\delta_{g+1}\}$ is a
$\delta$-sequence of type $A$ and that it is contained in
$S_{\nu,\infty}$. Now, we are going to prove that the semigroup
$S_{\nu,\infty}$ is generated by $\Delta$.
%$\{\delta_0,\delta_1,\ldots,\delta_g,\delta_{g+1}\}$.

Let us consider $f=f(x,y) \in k[x,y] \setminus \{0\}$ and set $f(x,y)=v^{-\deg (f)}
\bar{f} (u,v)$. To compute $\nu(f)$, one can use the procedure
described in the Remark 3.2 and write $f$ as an element of the
Laurent power series ring in the indeterminate $t$ with coefficients
in the field $k(s)$, $\mathfrak{L}_s(t)$. Indeed, by the proof of
Theorem 1 in \cite{belga}, it holds that
\[
\bar{f}(u(t,s),v(t,s))=(\mathfrak{q}_0 + \mathfrak{q}_1 p(s) + \cdots + \mathfrak{q}_r p(s)^r) t^{\nu(\bar{f})} + \mathfrak{r},
\]
where $\mathfrak{q}_i= \varrho_i s^{b_i}$, $b_i \geq 0$ and
$\varrho_i \in  k$, some $\mathfrak{q}_i$
does not vanish  and also the exponents of $t$ in $\mathfrak{r}
\in \mathfrak{L}_s(t)$ are always larger than  $\nu(\bar{f})$. On
the other hand, since $v$ is a curvette of any general element of
$\nu$,  $v(t,s) = \lambda_1 s^{a_1} t^{\nu(v)} + r_1$, $\lambda_1 \neq 0$, as we have
said in the Remark 3.2. As a consequence,  $f(x,y)$ can be written
as an element in $\mathfrak{L}_s (t)$ whose first jet in $t$ is
\begin{equation}
\label{aa}
(\mathfrak{m}_0 + \mathfrak{m}_1 p(s) + \cdots + \mathfrak{m}_r p(s)^r) t^{\nu(f)},
\end{equation}
$ \mathfrak{m}_i = s^{c_i} \mathfrak{q}_i$ for some $c_i \in \mathbb{Z}$, $0 \leq i \leq r$, which allows us to get the value $\nu(f)$.
An analogous situation happens for the polynomials $q_i=q_i (x,y)$, whose first jets are
\begin{equation}
\label{bb}
\gamma_i s^{d_i} t^{\nu(q_i)}, \; 0 \neq \gamma_i \in k, d_i \in \mathbb{Z},  \mbox{ for } q_i, \; 0 \leq i \leq g, \; \mbox{ and } s^d p(s) t^{\nu(q_{g+1})}, d \in \mathbb{Z}, \mbox{ for } q_{g+1}.
\end{equation}

Now, using \cite[Sect. 7]{pa2} and considering the class given by $f$ in $k[x,y]/(q_{g+1})$, we get that $f + (q_{g+1}) = \sum_{k=0}^d \xi_k \prod_{i=0}^g q_i^{s_{ik}} + (q_{g+1}) $, $ \xi_k \in k$ and the exponents satisfy the conditions $0 \leq s_{ik} < n_i$, ($1 \leq i \leq g$, $0 \leq k \leq d$). So,
 \begin{equation}
 \label{cc}
 f = \sum_{k=0}^d \left ( \xi_k \prod_{i=0}^g q_i^{s_{ik}} \right ) + q_{g+1}(x,y) f_1(x,y),
 \end{equation}
 where $f_1(x,y) \in k[x,y]$.

 Due to that the semigroup spanned by $\Gamma$ is telescopic \cite[Rem. 3.8]{far-cam}, the values $\nu(\prod_{i=0}^g q_i^{s_{ik}})$ are different, and this fact, together with the one expressed in (\ref{bb}), proves that $\nu(f)$ is the minimum of the values by $\nu$ of the summands involved in the expression (\ref{cc}).

 Finally, if $\nu(f) =  \nu(\prod_{i=0}^g q_i^{s_{ik}})$ for some
 $k$, $1 \leq k \leq d$, then the result is proved since
 $\nu(q_i) = - \delta_i$. If not, $\nu(f) = \nu(q_{g+1}(x,y) f_1(x,y))$ and
 we could repeat for $f_1$ the same procedure made with $f$ at most $r$
 times because otherwise using (\ref{bb}) and (\ref{aa}) we would obtain
 different expressions of $f$ in $\mathfrak{L}_s(t)$ which is not possible.
 Thus, it must happen that $\nu(f) = \nu (\prod_{i=0}^g q_i^{p_i} q_{g+1}^{p}) $
 for some non-negative integers $p_i$, $ 0 \leq i \leq g$ and  $p$, what concludes the proof in this case.

{\it Assume now that $\nu$ is a valuation of type} B. Let $C$ be the
plane curve having only one place at infinity whose successive
strict transforms (of the germ of $C$ at $p$) pass through all
points in ${\mathcal C}_{\nu}$.
%Making similar assumptions on $C$
As in the previous case, we can consider a family of approximates $\{
q_i(x,y)\}_{i=0}^g$ attached to the curve $C$ and an associated
$\delta$-sequence in $\mathbb{N}_{>0}$,
$\Gamma=\{\delta_0,\delta_1,\ldots,\delta_g\}$, such that
$S_{C,\infty}$ is generated by $\Gamma$ and the numbers $\delta_i$
satisfy the equalities given in Definition \ref{aprox}. Let $f(x,y)
\in k[x,y] \setminus \{0\}$ and suppose that $q(x,y) \in k[x,y]$ gives an equation
for $C$. With the notation as in (\ref{equality}), set
$f(x,y)=v^{-\deg (f)} \bar{f} (u,v)$ and analogously
$q(x,y)=v^{-\deg (q)} \bar{q} (u,v)$. Now, $\bar{f} (u,v) = \bar{q}
(u,v)^s \bar{l}(u,v)$, where $s \geq 0$, $\bar{l}(u,v) \in K[u,v]$
and $ \bar{q} (u,v)$ does not divide $\bar{l} (u,v)$.
%As above, one can see that there exists a polynomial of
%the form $\bar{m}(u,v) = u^{s_0} v^{s_1} \prod_{i=2}^g
%\overline{q}_i^{s_i} (u,v)$, $s_i \geq 0$ for any $i \geq 0$ such
%that $$ \nu_{C,p}(\bar{m}(u,v)) = \nu_{C,p}(\overline{f}(u,v)) - a
%\nu_{C,p} (\bar{q}(u,v))$$ and $\deg (\bar{m}(u,v)) \leq \deg
%(\bar{f}(u,v)) - a \deg (\bar{q}(u,v))$. Taking into account that
%for any element $h \in R$ such that $h (u,v) = \bar{q} (u,v)^a
%g(u,v)$, where $\bar{q} (u,v)$ does not divide $g(u,v)$, it holds
%that $\nu(h (u,v)) = (a, \nu_{C,p} (g(u,v))$, we get
%\[
%\nu(f(x,y))= \nu \left(x^{\deg (\bar{f})- a \deg (\bar{q})- \deg
%(\bar{m}) } q^a(x,y) m(x,y) \right)=
%\]
Recall that that for any element $h \in R$ such that $h (u,v) = \bar{q} (u,v)^a
 g(u,v)$, where $\bar{q} (u,v)$ does not divide $g(u,v)$, it holds
 that $\nu \left(h (u,v)\right) = \left(a, \nu_{C,p} (g(u,v) \right)$. Consider the polynomial $l(x,y) = v^{\deg(l)} \bar{l}(u,v)$, then  it happens that $\nu \left(l(x,y) \right) = \left(0, \nu_{C,p} (v^{\deg(l)} \bar{l}(u,v)) \right) =(0, \nu_{C,p} (l(x,y)))$. On the other hand, $\nu \left(q(x,y) \right)= \nu \left(v^{\deg(q)} \bar{q}(u,v) \right) =(0, - \delta_0^2) + (1,0)$, and this proves

\[
\nu \left(f(x,y) \right) = \nu \left(q^s (x,y) \right) + \nu \left( l(x,y) \right) = -s(-1,\delta_0^2) - \sum_{i=0}^g
s_i (0,\delta_i)
\]
for some $s_i \geq 0$. That is $\Delta := \{(0,\delta_0),
(0,\delta_1), \ldots, (0,\delta_g), (-1,\delta_0^2)\}$ is a
$\delta$-sequence of type B that spans $S_{\nu, \infty}$.

{\it The proof for valuations of type} C runs as follows. Set
$\{\nu_j\}_{j=1}^\infty $ a family of divisorial valuations that
define $\nu$. Consider $f(x,y) \in k[x,y]\setminus \{0\}$ and the corresponding
polynomial $\overline{f}(u,v)$ defined in (\ref{equality}). Pick a
valuation $\nu_{j_0}$ such that the intersection set between the
sequence of infinitely near points of any branch  of the germ
given by $\overline{f}(u,v)$ and the set $\mathcal{C}_\nu$
coincides with the one with $\mathcal{C}_{\nu_{j_0}}$. Let $C$ be
a curve having only one place at infinity and such that its germ
at $p$ gives a general element of $\nu_{j_0}$. Consider a family
of approximates $\{ q_i(x,y)\}_{i=0}^g$  for $C$. As above, we can suppose
%we can make assumptions on $C$ in such a way that
that $S_{C,\infty}$ is
spanned by a $\delta$-sequence in $\mathbb{N}_{>0}$,
$\{\delta_i\}_{i=0}^g$, which is related with the approximates as
in Definition \ref{aprox}. Remark \ref{nota} proves that
$\nu_{j_0} (f(x,y)) = \nu_{C,p} (f(x,y))$ and so, there exist
 positive integers $s_0,s_1,\ldots,s_g$ such that $\nu_{j_0} (f) =
\nu_{j_0} (\prod_{i=0}^g q_i^{s_i} (x,y))$. For analytically
irreducible elements $h$ in $R$, divisorial valuations $\nu_{j_0}$
satisfy a Noether formula, that is $\nu_{j_0} (h) = \sum_{j=0}^r
\nu_{j_0} (m_j) e(p_j)$, where $r$ gives the number of common
points $\{p_j\}_{j=0}^r$ between the sequences of infinitely near
points relative to $h$ and to $\nu$ and $e(p_j)$ is the
multiplicity of the strict transform of the germ given by $h$ at
the point $p_j$. The same Noether formula is also true for type C
valuations with the values described for $\nu$ in Section
\ref{ronaldo} (see \cite[Sect. 3]{galmon}) and this proves that
$\nu (f) = \nu (\prod_{i=0}^g q_i^{s_i} (x,y))$. Therefore
$S_{\nu,\infty}$ is spanned by the values $-\nu(q_i(x,y))\in
\mathbb{Z}^2$, $0\leq i\leq g$, which constitute a
$\delta$-sequence $\Delta$ of type C that can be computed from the
above $\delta$-sequence in $\mathbb{N}_{>0}$ as is described in
Definition \ref{buena}.

{\it Suppose now that $\nu$ is a valuation of type} D. Let
$\{\nu_j\}_{j=1}^\infty$ be a family of divisorial valuations at
infinity that define $\nu$. Without loss of generality, we can
assume that, for all $j$, $\nu_j$ is the divisorial valuation
defined by the exceptional divisor associated with a blowing-up
centered at a satellite point and, moreover, all the free points
in ${\mathcal C}_{\nu}$ also belong to ${\mathcal C}_{\nu_j}$. Set
$\{\delta_i^j\}_{i=0}^g$  a normalized $\delta$-sequence in
$\mathbb{N}_{>0}$ associated with a curve having only one place at
infinity which provides a general element of $\nu_j$.  As we have
said, for $f(x,y) \in k[x,y] \setminus \{0\} $ and if $j$ is large enough, the
value of $f$ by the valuation $\nu_j$  (normalized as we explained
at the end of Section \ref{ronaldo}) is $\nu_j (f) = \sum_{i=0}^g
s_i \delta_i^j$, $s_0,s_1,\ldots,s_g$ being non-negative integers.
Thus,
$$\nu(f(x,y)) = \lim_{j \rightarrow
\infty} \nu_j (f) =\lim_{j \rightarrow \infty}  \sum _{i=0}^g s_i
\delta_i^j.
$$

For simplicity assume that $g\geq 2$ (all works similarly when
$g=1$). By Corollary \ref{coro},  one has that when $i<g$,
$\delta_i^j=\delta_i^{j'}$ for whichever pair $j,j'$ of indices of
the divisorial valuations. Moreover,
from the  paragraph in Section \ref{ronaldo} corresponding to this
type of valuations, it is straightforward to deduce that
$\delta_g:=\lim_{j \rightarrow \infty} \delta_g^j$ is a
non-rational number. Therefore, $S_{\nu,\infty}$ is spanned by the
$\delta$-sequence of type D given by $\Delta:=
\{\delta_0^j,\delta_1^j,\ldots,\delta_{g-1}^j,\delta_g\}$.

Finally, {\it whenever $\nu$ is a valuation of type} E, taking
$f(x,y)$, $\{\nu_j\}_{j=1}^\infty $ and $\nu_{j_0}$ as in the case
of valuations of type C, it holds that $\nu(f) = \nu_{j_0} (f)$
(here, as in the above paragraph, the valuations $\nu_{j}$ are
also normalized). From this fact and arguing in a similar way as
above it follows immediately that $S_{\nu, \infty}$ is spanned by
a $\delta$-sequence $\Delta$ of type E, whose first elements
coincide with the normalized $\delta$-sequence of the
$\delta$-sequence in $\mathbb{N}_{>0}$ defined by a curve with
only one place at infinity whose germ in $p$ gives a general
element of $\nu_{j_0}$.\findemo

\begin{rem}
\label{r42}
{\rm

Suppose that ${\rm char} (k)=0$ and let $\nu$ be a valuation at
infinity of type A (notice that $\nu$ satisfies the hypothesis of
Theorem \ref{UNICO}). In this case the sign of $\delta_{g+1}$ has a
nice geometrical interpretation. Indeed, set
$\Delta=\{\delta_0,\delta_1,\ldots,\delta_g,\delta_{g+1}\}$ a
$\delta$-sequence of type A generating the semigroup
$S_{\nu,\infty}$ and let $f(x,y)=0$ be an equation defining a curve
$C$ having only one place at infinity whose germ at $p$ is a general
element for $\nu$. In addition, assume that $f(x,y)$ and
$\{\delta_0,\delta_1,\ldots,\delta_g\}$ are as in Definition
\ref{aprox} (in particular, $\delta_0$ is the degree of $C$). Let us
consider the infinite sequence
$$\cdots \rightarrow X_i\rightarrow X_{i-1}\rightarrow \cdots
\rightarrow X_1\rightarrow X_0:=\gp^2$$ such that $X_1\rightarrow
\gp^2$ is the blowing-up of $\gp^2$ at $p_0:=p$ and, for each $i\geq
2$, $X_{i}\rightarrow X_{i-1}$ is the blowing-up at the intersection
point $p_{i-1}$ of the strict transform of $C$ in $X_{i-1}$ and the
exceptional divisor of the preceding blowing-up. Assume that the
first $r+1$ blowing-ups induce the minimal embedded resolution of
the singularity of $C$ at $p$ (as in (\ref{infiniteseq})). Let
$\Lambda$ be the pencil of plane curves of $\gp^2$ whose equations
in the affine chart $Z\not=0$ are $f(x,y)=\lambda$, where
$\lambda\in k$. By a result of Moh \cite{moh}, all the curves in
$\Lambda$ have only one place at infinity and there exists an integer $\mathbf{m} \geq r$ such
that their strict transforms pass through $p_i$ for all
$i\leq \mathbf{m}$, $\mathbf{m}$ being the least integer such that the
composition of the first $\mathbf{m}+1$ blowing-ups eliminates the base
points of the pencil $\Lambda$. Thus
$\delta_0^2=\sum_{i=0}^{\mathbf{m}} m_{p_i}(C)^2$, where $m_{p_i}(C)$ denotes
the multiplicity at $p_i$ of the strict transform of $C$ in $X_i$.
By the proof of Theorem \ref{UNICO}, it happens
$$
\delta_{g+1}=-\nu({q}_{g+1}(x,y))=\delta_0^2-\nu(\bar{q}_{g+1}(u,v))=\delta_0^2-\sum_{i=0}^{\mathbf{n}} m_{p_i}(C)^2,
$$
where ${\mathcal
C}_{\nu}=\{p_i\}_{i=0}^\mathbf{n}$.

Therefore, if $\mathbf{n}\leq \mathbf{m}$
(respectively, $\mathbf{n}>\mathbf{m}$), then $\delta_{g+1}\geq 0$ (respectively, $<0$) and so
$S_{\nu,\infty}$ is (respectively, is not) a well-ordered semigroup (with
the natural ordering).

}
\end{rem}

\begin{rem}
{\rm Theorem \ref{UNICO} holds under certain conditions on the
characteristic of the field $k$. When the plane valuation at
infinity $\nu$ is either of type A or B, it has no difficulty to
check that condition. Otherwise and when char$(k) \neq 0$, it
cannot be easy to decide whether the cited condition happens.
Notwithstanding, it is not difficult to find examples where to
check that condition is simple. For instance, following
\cite[Sect. 4.3.4]{galmon}, a valuation of type $E$ can be
obtained considering, as a $\delta$-sequence, the limit of the
normalized $\delta$-sequences of the $\delta$-sequences in
$\mathbb{N}_{>0}$: $\{5,3\}$, $\{5\cdot 2,3 \cdot 2, 3 \cdot 3
\}$, $\{5\cdot 2^2,3 \cdot 2^2, 3^2 \cdot 2, 3^3 \}$, \ldots (we
have used $\{5,3\}$ to start and the value $z=2$ with the notation
in \cite[Sect. 4.3.4]{galmon}). Then, Theorem \ref{UNICO} happens
when the characteristic of the field is different from 2. Let us
see another example corresponding to a type C valuation. Let $\nu$
be the valuation of type C whose $\delta$-sequence is
$\{(2,1),(1,0) \}$. Then the sequence of values
$(-\nu_i(x),-\nu_i(y))$ described in the Theorem \ref{UNICO} is an
infinite subset of the set $\{(2j+1,j)\}_{j \geq 2}$. Since
$\gcd(2j+1,j)=1$ for any index $j$, our result happens for any
characteristic of the field $k$. }\end{rem}

\begin{rem}
{\rm The converse of Abhyankar-Moh Theorem for plane valuations at
infinity, without any restriction on the characteristic of the
field, is true. It is proved for valuations of types C, D and E in
\cite[Th. 4.9]{galmon}. Let us see it for the types A and B.
Assume that $\Delta$ is a $\delta$-sequence either of type A or B.
Set $\Delta'=\{\delta_0,\delta_1,\ldots,\delta_g\}$ the
$\delta$-sequence in $\mathbb{N}_{>0}$ that allows us to provide
$\Delta$. Consider the unique expression
\[
\label{siete}
n_i \delta_i = \sum_{j=0}^{i-1} a_{ij} \delta_j,
\]
where $a_{i0} \geq 0$ and $0 \leq a_{ij} < n_j$, for $1 \leq j \leq i-1$. Let us define the following polynomials: $q_0:=x$, $q_1:=y$ and, for $1 \leq i
\leq g$,
\[
 \label{misq}
q_{i+1}:= q_{i}^{n_{i}} - t_i \left ( \prod_{j=0}^{i-1} q_j^{a_{ij}} \right),
\]
where $t_i \in k \setminus \{0\}$ are arbitrary. The curve
$C_{\Delta'}$ defined by $q_{g+1}(x,y)=0$ has only one place at
infinity, the set $\{q_i\}_{i=0}^g$ is a family of approximates for
it and $S_{C_{\Delta'},\infty}$ is generated by $\Delta'$ (see
\cite{galmon} and references therein for more details).

Then the valuation $\nu$ defined by the  set
$\{p_i\}_{i=1}^\infty$ of infinitely near points attached to the
resolution of $C_{\Delta'}$ is a valuation at infinity of type B
whose semigroup at infinity is spanned by the $\delta$-sequence
$\Delta$ of type B. In addition, if we take the valuation $\nu$
defined by the finite subset $\{p_i\}_{i=1}^n$ of the above set
$\{p_i\}_{i=1}^\infty$ such that $n=n_g \delta_g - \delta_{g+1} +
i_0$, $\delta_{g+1}$ being the last element in $\Delta$ and $i_0$
the largest index such that $p_{i_0}$ is a satellite point, then
the semigroup at infinity of $\nu$ is spanned by the $\delta$-sequence
$\Delta$ of type A.}
\end{rem}

%\newpage

\end{document}